\documentclass[12pt]{article}
\usepackage{latexsym}
\usepackage{amsmath}
\usepackage{amssymb}
\usepackage{amsfonts}
\usepackage{graphics,graphicx}

\usepackage{url}

\def\claim#1{\begin{trivlist}\item[\hskip\labelsep\bf#1]\it}
\def\endclaim{\end{trivlist}}

\numberwithin{equation}{section}

\newtheorem{theorem}{Theorem}[section]
\newtheorem{lemma}[theorem]{Lemma}

\newcommand{\eproof}{{\mbox{\ }~\hfill
\mbox{\large $\Box$} \par \vskip 10pt}}
\newcommand{\pf}{\noindent{\bf Proof}}

\title{ Unique continuation property for anomalous slow diffusion equation}
\author{Ching-Lung Lin\thanks{Department of Mathematics, NCTS, National Cheng-
Kung University, Tainan 701, Taiwan.  Partially supported by the
National Science Council of Taiwan. (Email:
cllin2@mail.ncku.edu.tw)}\qquad Gen Nakamura\thanks{Department of
Mathematics, Inha University, Incheon 402-751, South Korea.
Partially supported by Korea Research Foundation.
(Email: 213028@inha.ac.kr)}}

\date{}

\begin{document}
\renewcommand{\theequation}{\thesection.\arabic{equation}}
 \maketitle
\begin{abstract}
A Carleman estimate and the unique continuation of solutions for an anomalous diffusion equation with fractional time derivative of order $0<\alpha<1$ are given. The estimate is derived via some subelliptic estimate for an operator associated to the anomalous diffusion equation using calculus of pseudo-differential operators.

\end{abstract}
\section{Introduction}\label{sec1}
In this paper we are concern with the UCP (unique continuation property) of solutions of anomalosu diffusion equation. Even in the case discussing UCP of solution in $H^{\alpha, 2}(\Omega\times (0,T))$ (see below in this section for its definition) giving zero Cauchy data on a small part $\Gamma$ of the $C^2$ boundary $\partial\Omega$ of a domain $\Omega\subset{\Bbb R}^n$ over some time interval, we can always consider the $0$ extension of the solution outside $\Omega$ in a neighborhood of $\Gamma$. Hence assuming $0\in\partial\Omega$ without loss of generality we only need to consider the following for UCP of solutions of anomalous diffusion equation. That is let $\hat{y}=(\hat{y}_1,\cdots,\hat{y}_{n-1},0)\in \mathbf{R}^n$ and $\omega=\{(y_1,\cdots,y_n):\hat{y}_j-l<y_j<\hat{y}_j+l,-l<y_n\leq0, 1\leq j\leq n-1\}\subset\mathbf{R}^n$ with $l>0$, consider the following equation for $0<\alpha<1$
\begin{equation}\label{1.1}
\begin{cases}
\begin{array}{l}
\partial^{\alpha}_{t}u(t,y)-\Delta_y u(t,y)
=l_1(t,y;\nabla_y)u(t,y),\\
u(t,y)=0\quad (t\leq0),\\
u(t,y)=0 \quad (y\in \omega, 0<t<T),
\end{array}
\end{cases}
\end{equation}
where $l_1$ is a linear differential operator of order $1$ and $\partial^{\alpha}_{t}u$ is the fractional derivative of $u$
in the Caputo sense which is defined by
$$\partial^{\alpha}_{t}u(t,y)=\frac{1}{\Gamma(1-\alpha)}((t^{-\alpha}H(t)\otimes\delta_y)\ast\partial_tu)(t,y),$$
where $H(t)$ is the Heaviside function and $\delta_y$ is the Dirac function with singularity at $y$.

\medskip

The anomolous diffusion equation was first studied in material science and the exponent $\alpha$ of $\partial^{\alpha}_{t}u$ in this equation is an index which describes the long time behavior of the mean square displacement $<x^2(t)>\sim\mbox{\it positive const.}\,t^\alpha$ of a diffusive particle $x(t)$ describing anomalous diffusion on fractals such as some amorphous semiconductors or strongly
porous materials (see \cite{Anh}, \cite{Metzler} and references therein).

Recently a strong inertia to the study of anomolous diffusion equation came from a study in environmental science. It showed by an experiment that the spread of pollution in soils cannot be modeled correctly by the usual diffusion equation, but it can be modeled by
an anomalous diffusion equation (see \cite{hat1}, \cite{hat2}). The Cauchy problem and intial boundary value problem for the anomalous diffusion equation have been studied by many people (see \cite{A}, \cite{E} and the references therein).

The aim of this paper is to give a Carleman estimate of solutions of anomalous diffusion equation which enables us to have  UCP of its solutions for any $\alpha\,(0<\alpha<1)$ and $n\in{\Bbb N}$. UCP is a key to the study of control problem and inverse problem for this equation. It can give
the approximate boundary controlability for the control problem and it is very important for inverse problem if one wants to develop for instance linear sampling type reconstruction scheme to identify unknown objects such as cracks, cavities and inclusions inside an anomalous diffusive medium. Some Carleman estimates have been given for some special cases. That is for
$\alpha=1/2$, a Carleman estimate was given in \cite{XCY}, \cite{ZX}, \cite{ZY} for $n=1$ and \cite{CLN} for $n=2$ via that for the operator $\partial_t+\Delta^2$ with some
lower order terms.

Looking at the symbol of the anomalous diffusion equation, we can say that it is basically semi-elliptic. Based on this, we adapt Treve's argument (\cite{Treve}) to derive a Carleman estimate via some subelliptic estimate to the equation conjugated by $e^{\tau_0 t}$ with $\tau_0$ and transformed by a Holmgren transform. Our method has a potential to be applied to space time fractional diffusion equations (\cite{Hanyga}) and some fractional derivative visco-elastic equations (\cite{Lu}).

\par In order to state our results, let $\mathcal{S}(\mathbf{R}_t^1\times\mathbf{R}_x^n)$ and $\mathcal{S}'(\mathbf{R}_t^1\times\mathbf{R}_x^n)$ be
the set of rapidly decreasing functions in $\mathbf{R}_t^1\times\mathbf{R}_x^n$ and its dual space, respectively.
Then, for $m,s\in \mathbf{R}$, $v=v(t,x)\in \mathcal{S}'(\mathbf{R}_t^1\times\mathbf{R}_x^n)$, belongs to the function space
$H^{m,s}(\mathbf{R}_t\times\mathbf{R}_x^n)$ if
$$\|v\|^2_{H^{m,s}}:=\iint(1+|\xi|^s+|\tau|^m)^2|\hat{v}|^2d\tau d\xi$$
is finite, and $\|v\|^2_{H^{m,s}}$ denotes the norm of $v$ of this function space, where
$\hat{v}$ is the Fourier transform defined by
$$\mathcal{F}(v)(\tau,\xi)=\hat{v}(\tau,\xi)=\iint e^{-it\tau-ix\cdot\xi}v(t,x)dxdt.$$
Further, we define $H^{m,s}(\Omega\times(0,T))$ as the restriction of $H^{m,s}(\mathbf{R}_t\times\mathbf{R}_x^n)$ to $\Omega\times(0,T)$.

In this paper, we will give a Carleman estimate and the following UCP of solutions of \eqref{1.1}.

\begin{theorem}\label{thm1.1}
Let $u\in H^{\alpha,2}(\mathbf{R}^{1+n})$ satisfy \eqref{1.1}. Then
$u$ will be zero across $y_n=0$.
\end{theorem}

The rest of this paper is organized as follows. We compute the principal symbol of the anomalous diffusion operator which undergone Holmgren type transformation and then multiplied by an exponential function of time variable, and analyze its properties in Section 2. In Section 3, we derive some subelliptic estimate for some pseudo-differential operator associated to this operator. Then by using this subelliptic estimate, the Carleman estimate is derived in
Section 4. Finally in the last section, we give UCP of solutions to the anomalous diffusion operator.

\section{Change of variables and principal symbol}\label{sec2}

UCP is a local property of solutions of \eqref{1.1}. Hence we consider $\theta(t)\kappa(y_n)u(t,y)$ with
$\kappa(y_n)\in C^\infty(\mathbf{R})$ and $\theta(t)\in C^\infty(\mathbf{R})$ defined by
\begin{equation*}
\begin{array}{l}
\kappa(y_n)=
\begin{cases}
0,\quad  y_n\leq -2l/3,\\
1, \quad y_n\geq -l/3
\end{cases}
\end{array}
\end{equation*}
and
\begin{equation*}
\begin{array}{l}
\theta(t)=
\begin{cases}
1, \quad t\leq T-\epsilon,\\
0, \quad t\geq T-\epsilon/2
\end{cases}
\end{array}
\end{equation*}
for small positive $\epsilon$. By abusing the notation
$u(t,y)$ to denote $\theta(t)\kappa(y_n)u(t,y)$, this $u(t,y)$ satisfies
\begin{equation}\label{2.1}
\begin{cases}
\begin{array}{l}
\partial^{\alpha}_{t}u(t,y)-\Delta_y u(t,y)
=l_1(t,y;\nabla_y)u(t,y),\\
u(t,y)=0\quad (t\leq0),\\
u(t,y)=0 \quad (y_n\leq 0\quad {\rm or}\quad t\geq T).
\end{array}
\end{cases}
\end{equation}
Then our aim is to show that $u$
considered on $\{y_n\geq 0\}\cap \{t\leq T-\epsilon\}$ where we have $\theta(t)\kappa(y_n)=1$
will be zero across $y_n=0$.

To have UCP of solutions of \eqref{2.1} in the $y_n$-direction by a Carleman estimate,
we use the change of variables $x'=y'-\hat{y}', x_n=y_n+|y'-\hat{y}'|^2+\frac{X}{T}(t-T), t=t$, where $x'=(x_1,\cdots,x_{n-1})$, $y'=(y_1,\cdots,y_{n-1})$, $\hat{y}'=(\hat{y}_1,\cdots,\hat{y}_{n-1})$ and $X$ is a small positive constant which will be determined later.
This is a Holmgren type transformation.
By $\partial_t=\frac{X}{T}\partial_{x_n}+\partial_t$,
$ \partial_{y_j}=\partial_{x_j}+2x_j\partial_{x_n}$ for $j=1,\cdots,n-1$, we have
\begin{equation}\label{2.2}
\begin{array}{l}
\partial_{t}^\alpha u(t,y)=\frac{1}{\Gamma(1-\alpha)}\int_0^t(t-\eta)^{-\alpha}\partial_\eta u(\eta,y)d\eta\\
=\frac{1}{\Gamma(1-\alpha)}\int_0^t(t-\eta)^{-\alpha}\partial_\eta u(\eta,x)d\eta
+\frac{X}{T\Gamma(1-\alpha)}\int_0^t(t-\eta)^{-\alpha}\partial_{x_n}u(\eta,x)d\eta,
\end{array}
\end{equation}
where $\Gamma(\cdot)$ denotes the Gamma function and $u(t,x)$ with $t=\eta$ on the right hand side of \eqref{2.2}
is the push forward of $u(t,y)$ by the change of variables. We further abuse the notation $u(t,x)$ to denote this $u(t,x)$ multiplied by
$e^{\tau_0 t}$ with $\tau_0<0$. Then, this new $u(t,x)$ also satisfies
$${\rm supp}\, u\subset\{x_n\geq -X\}.$$

For further arguments, we define a function space $\dot{H}^m_{\alpha}(\overline{\mathbf{R}_+^{1+n}})$ for $m\in \mathbf{R}$ as follows.
First, let $(t,x)\in\mathbf{R}^{1+m}$ and denote the open half of $\mathbf{R}^{1+m}$ space defined by $t>0$ and complement of its closure $\bar{\mathbf{R}}_+^{1+m}$
by $\mathbf{R}_+^{1+m}$ and $\mathbf{R}_-^{1+m}$, respectively. If $E$ is a space of distributions in $\mathbf{R}^{1+m}$ we
use the notation $\bar{E}(\mathbf{R}_+^{1+m})$ for the space of restrictions to $\mathbf{R}_+^{1+m}$ of elements in $E$ and we write
$\dot{E}(\bar{\mathbf{R}}_+^{1+m})$ for the set of distributions in $E$ supported by $\bar{\mathbf{R}}_+^{1+m}$.

Next, set
$$\Lambda_{\alpha}^m(\tau,\xi)=((1+|\xi|^2)^{1/\alpha}+i\tau)^{m\alpha/2} \quad \mbox{for}\,\,m\in \mathbf{R}$$
and define a pseudo-differential operator $\Lambda_{\alpha}^m(D_t,D_x)$ by $\Lambda_{\alpha}^m(D_t,D_x)\phi:=\mathcal{F}^{-1}(\Lambda_{\alpha}^m(\tau,\xi)\hat{\phi})$, $\phi\in \dot{\mathcal{S}}(\overline{\mathbf{R}_+^{1+n}})$, where
$\mathcal{F}^{-1}$ is the inverse Fourier transform. Then
$\Lambda_{\alpha}^m(D_t,D_x)$ is a continuous map in $\dot{\mathcal{S}}(\overline{\mathbf{R}_+^{1+n}})$
with inverse $\Lambda_{\alpha}^{-m}(D_t,D_x)$. (See Theorem B.2.4 in \cite{Hormander}.) Based on this we define $\dot{H}^m_{\alpha}(\overline{\mathbf{R}_+^{1+n}})$ by $\dot{H}^m_{\alpha}(\overline{\mathbf{R}_+^{1+n}}):=\Lambda_{\alpha}^{-m}L^2(\mathbf{R}_+^{1+n})$,
where we considered $L^2(\mathbf{R}_+^{1+n})$ as a subset of $\dot{\mathcal{S}}'(\overline{\mathbf{R}_+^{1+n}})$.

Now denote by $P(t,x;D_t,D_x)$ the operator $e^{\tau_0 t}(\partial_t^\alpha-\Delta_y)$ in terms of the coordinates $(t,x)$ and its total symbol by
$p(t,x;\tau,\xi)$, where $D_t=-\sqrt{-1}\partial_t$, $D_x=-\sqrt{-1}\partial_x$. Then it is easy to see
\begin{equation*}
\begin{array}{rl}
p(t,x;\tau,\xi)&=(i(\tau+i\tau_0))^\alpha+\Sigma_{j=1}^{n-1}(\xi_j+2x_j\xi_n)^2+\xi_n^2+\frac{X}{T}i^\alpha(\tau+i\tau_0)^{\alpha-1}\xi_n\\
&=(i(\tau+i\tau_0))^\alpha+|\xi'|^2+4g\xi_n+f\xi_n^2+\frac{X}{T}i^\alpha(\tau+i\tau_0)^{\alpha-1}\xi_n,
\end{array}
\end{equation*}
where  $g=g(x',\xi')=\Sigma_{j=1}^{n-1}x_j\xi_j$ and $f=f(x')=1+4|x'|^2$.
Also, by the Payley-Wiener type theorem, $P(t,x;D_t,D_x)$ maps $\dot{H}^m_{\alpha}(\overline{\mathbf{R}_+^{1+n}})$
into $\dot{H}^{m-2}_{\alpha}(\overline{\mathbf{R}_+^{1+n}})$ for any $m\in \mathbf{R}$.

We want to derive a Carleman estimate for the operator $P(t,x;D_t,D_x)$ using some subelliptic estimate based on the idea by F. Treve (\cite{Treve}). For that let $\psi=\frac{1}{2}(x_n-X)^2$ and consider the symbol $p(x;\tau,\xi+i|\sigma|\nabla\psi)$ over $\mathbf{R}^{n+1}\times\mathbf{R}_z$ to define a pseudo-differential operator $P_\psi(x,D_t,D_{x},D_z)$ by
$$P_\psi=P_\psi(x,D_t,D_{x},D_z)=p(x,D_{t},D_x+i|D_z|\nabla \psi)$$
which can be given for any compactly supported distribution $v$ in $\mathbf{R}^{n+1}\times\mathbf{R}_z$ by
\begin{equation}\label{2.3}
\begin{array}{rl}
P_\psi v(z,t,x)=\int e^{i(x\cdot\xi+t\tau+z\sigma)}p(x;\tau,\xi+i|\sigma|\nabla\psi)\hat{v}(\sigma,\xi,\tau)d\sigma d\tau d\xi.
\end{array}
\end{equation}
We denote the principal symbol of $P_\psi$ by $\tilde{p}_\psi$ which is given by
\begin{equation*}
\begin{array}{rl}
\tilde{p}_\psi=&(i(\tau+i\tau_0))^\alpha+|\xi'|^2+4g\xi_n+f\xi_n^2-f|\sigma|^2(x_n-X)^2\\
&+i4g(x_n-X)|\sigma|+i2f\xi_n(x_n-X)|\sigma|.
\end{array}
\end{equation*}
It is important to note here that $\tilde{p}_\psi$ is independent of $t$ and $z$. This gives the advantage in computing the Poisson bracket of $\tilde{p}_\psi$.

By the definition of Poisson bracket
\begin{equation}\label{Poisson bracket}
\{\Re \tilde{p}_\psi, \Im \tilde{p}_\psi\}=\Sigma_{j=1}^{n}(\partial_{\xi_j}\Re \tilde{p}\cdot\partial_{x_j}\Im \tilde{p}
-\partial_{x_j}\Re \tilde{p}\cdot\partial_{\xi_j}\Im \tilde{p})
\end{equation}
with the real part $\Re \tilde{p}_\psi$ and imaginary part $\Im \tilde{p}_\psi$ of $\tilde{p}_\psi$ given by
\begin{equation*}
\begin{array}{ll}
\Re \tilde{p}_\psi=\Re(i(\tau+i\tau_0))^\alpha+|\xi'|^2+4g\xi_n+f\xi_n^2-f(x_n-X)^2|\sigma|^2\\
\Im \tilde{p}_\psi=\Im(i(\tau+i\tau_0))^\alpha+4g(x_n-X)|\sigma|+2f\xi_n(x_n-X)|\sigma|.
\end{array}
\end{equation*}
Note that there are no terms related with $t$ and $z$ derivatives in \eqref{Poisson bracket}.
A direct computation gives
\begin{equation*}
\begin{array}{l}
\nabla_\xi\Re \tilde{p}_\psi=2(\xi',0)+4\xi_n(x',0)+4g(0',1)+2f\xi_n(0',1)\\
\nabla_x\Im \tilde{p}_\psi=4(x_n-X)|\sigma|(\xi',0)+4g|\sigma|(0',1)+16\xi_n(x_n-X)|\sigma|(x',0)+2f\xi_n|\sigma|(0',1)\\
\nabla_x\Re \tilde{p}_\psi=4\xi_n(\xi',0)+8\xi_n^2(x',0)-8(x_n-X)^2|\sigma|^2(x',0)-2f(x_n-X)|\sigma|^2(0',1)\\
\nabla_\xi\Im \tilde{p}_\psi=4(x_n-X)|\sigma|(x',0)+2f(x_n-X)|\sigma|(0',1),
\end{array}
\end{equation*}
where $0'$ denotes $\xi'=0$.
Thus
\begin{equation*}
\begin{array}{rl}
&\Sigma_{j=1}^{n}(\partial_{\xi_j}\Re \tilde{p}\cdot\partial_{x_j}\Im \tilde{p})\\
=&8(x_n-X)|\sigma||\xi'|^2+48g\xi_n(x_n-X)|\sigma|+64|x'|^2\xi_n^2(x_n-X)|\sigma|\\
&+16g^2|\sigma|+16fg\xi_n|\sigma|+4f^2\xi_n^2|\sigma|
\end{array}
\end{equation*}
and
\begin{equation*}
\begin{array}{rl}
&\Sigma_{j=1}^{n}(\partial_{x_j}\Re \tilde{p}\cdot\partial_{\xi_j}\Im \tilde{p})\\
=&16g\xi_n(x_n-X)|\sigma|+32|x'|^2\xi_n^2(x_n-X)|\sigma|-32|x'|^2(x_n-X)^3|\sigma|^3-4f^2(x_n-X)^2|\sigma|^3.
\end{array}
\end{equation*}
Observe that
$$
(i(\tau+i\tau_0))^\alpha=|\tau+i\tau_0|^\alpha[\cos(\alpha \arg(i(\tau+i\tau_0)))+i\sin((\alpha \arg(i(\tau+i\tau_0)))],
$$
we have for some $\epsilon_0>0$
\begin{equation}\label{3.4}
\begin{array}{l}
\Re (i(\tau+i\tau_0))^\alpha\geq\epsilon_0|\tau|^\alpha.
\end{array}
\end{equation}
From $\Re \tilde{p}_\psi=0$, $|x'|^2\leq X/4$ with $X<<1$ and \eqref{3.4}, we have
\begin{equation}\label{2.4}
\begin{array}{l}
(x_n-X)^2|\sigma|^2\gtrsim |\tau|^\alpha+|\xi|^2+\sigma^2.
\end{array}
\end{equation}
Also from \eqref{2.4}, we have
\begin{equation}\label{2.5}
\begin{array}{l}
\{\Re \tilde{p}_\psi, \Im \tilde{p}_\psi\}\gtrsim (|\tau|^\alpha+|\xi|^2+\sigma^2)^{3/2}.
\end{array}
\end{equation}

\section{Subelliptic estimates}\label{sec3}
In this section we will show the following subelliptic estimate for the operator $P_\psi$.
\begin{lemma}\label{lem3.2}
There exists a sufficiently small constant $z_0$ such that for all $u(t,x,z)\in C_0^{\infty}(U\times[-z_0,z_0])\cap\dot{\mathcal{S}}(\bar{\mathbf{R}}_+^{1+n+1})$, we have that
\begin{equation}\label{3.1}
\begin{array}{l}
\Sigma_{k+s<2}||h(D_z)^{2-k-s}\Lambda_\alpha^{s}D_z^k u||
\lesssim||P_\psi u||,
\end{array}
\end{equation}
where $h(D_z)=(1+D_z^2)^{1/4}$ and $U$ is a small open neighborhood of the origin in ${\Bbb R}^{1+n}$.
\end{lemma}
\pf.

Let us suppose there is a finite open covering of
$$|\xi|^2+\sigma^2+|\tau|^\alpha=1$$
and a subordinate smooth partition of unity $\chi_\nu(\xi,\tau,\sigma)$ homogeneous of degree $0$ in terms of the scaling
$$
(\xi,\tau,\sigma)\mapsto (\eta\xi,\eta^{2/\alpha}\tau,\eta\sigma)
$$
for $\eta>0$ such that
$\Sigma\chi_\nu^2=1$ and
\begin{equation}\label{3.2}
\begin{array}{l}
||\chi_\nu(D_x,D_t,D_z)h(D_z) u||^2_{H^{\alpha/2,1}}
\lesssim||P_\psi\chi_\nu u||_{L^2}^2+||u||^2_{H^{\alpha/2,1}}.
\end{array}
\end{equation}
Here we have abused the notation $\Vert u\Vert_{H^{m,s}}$ with $m,\,s\in\mathbf{R}$ to denote the norm for
$u(t,x,z)\in\dot{\mathcal{S}}(\bar{\mathbf{R}}_+^{1+n+1})$ given by
$$
\Vert u\Vert_{H^{m,s}}^2:=\int\int\int(1+|\tau|^m+|\xi|^s+|\sigma|^s)^2|\hat u|^2\,d\tau\,d\xi\,d\sigma.
$$
Summing up all $\nu$ on \eqref{3.2}, we can have from elementary estimates that
\begin{equation}\label{3.3}
\begin{array}{l}
||h(D_z) u||^2_{H^{\alpha/2,1}}
\lesssim||P_\psi u||_{L^2}^2
\end{array}
\end{equation}
which implies \eqref{3.1} for a small enough $z_0$.

So, it suffices to work microlocally. That is, we need to establish \eqref{3.2}.
Let's divide into two cases to show \eqref{3.2}. If $\chi_\nu$ is supported in a small neighborhood of $\sigma=0$,
then $|\xi|^2+|\tau|^\alpha\geq \delta_0>0$.
We recall that the principal symbol $\tilde{p}_\psi$ of $P_\psi$ is
\begin{equation*}
\begin{array}{rl}
\tilde{p}_\psi=&(i(\tau+i\tau_0))^\alpha+|\xi'|^2+4g\xi_n+f\xi_n^2-f|\sigma|^2(x_n-X)^2\\
&+i4g(x_n-X)|\sigma|+i2f\xi_n(x_n-X)|\sigma|.
\end{array}
\end{equation*}
A direct computation gives for small enough $X$ that
\begin{equation}\label{3.5}
\begin{array}{rl}
|\tilde{p}_\psi|\geq|\Re\tilde{p}_\psi|&=\Re (i(\tau+i\tau_0))^\alpha+|\xi'|^2+4g\xi_n+f\xi_n^2-f|\sigma|^2(x_n-X)^2\\
&\geq \epsilon_0(|\xi|^2+|\tau|^\alpha)-f|\sigma|^2(x_n-X)^2\\
&\gtrsim \delta_0\gtrsim |\xi|^2+\sigma^2+|\tau|^\alpha,
\end{array}
\end{equation}
where $\epsilon_0$ is the constant in \eqref{3.4} and $\delta_0$ is a small positive constant.
So, we can have a better result than \eqref{3.2} in this case.

On the other hand, if the support of $\chi_\nu$ is bounded away from $\sigma=0$,
then $\sigma^2\geq \delta_1(|\xi|^2+|\tau|^\alpha)$ with a positive constant $\delta_1$.
We write $||P_\psi\chi_\nu u||_{L^2}^2=(P_\psi\chi_\nu u,P_\psi\chi_\nu u)$ as the following.
\begin{equation}\label{3.6}
\begin{array}{rl}
||P_\psi\chi_\nu u||_{L^2}^2&=(P_\psi\chi_\nu u,P_\psi\chi_\nu u)\\
&=(P^*_\psi P_\psi\chi_\nu u,\chi_\nu u)\\
&=(P_\psi P^*_\psi\chi_\nu u,\chi_\nu u)+([P^*_\psi, P_\psi]\chi_\nu u,\chi_\nu u)\\
&=\big((I-\eta B)P_\psi P^*_\psi\chi_\nu u,\chi_\nu u\big)+\big(([P^*_\psi, P_\psi]+\eta BP_\psi P^*_\psi)\chi_\nu u,\chi_\nu u\big),
\end{array}
\end{equation}
where the principal symbol of the commutator $[P^*_\psi, P_\psi]$
is $[\overline{\tilde{p}_\psi}, \tilde{p}_\psi]=2\{\Re \tilde{p}_\psi, \Im \tilde{p}_\psi\}$ which has been already studied in Section \ref{sec2}, $\eta$ a large positive constant and $B=\Lambda^{-1/2}(\Lambda^{-1/2})^*$ with an elliptic pseudo-differential operator $\Lambda$ whose principal symbol is $(|\tau|^\alpha+|\xi|^2+\sigma^2)^{1/2}$.
From \eqref{2.5}, we have that the principal symbol of $[P^*_\psi, P_\psi]+\eta BP_\psi P^*_\psi$ satisfies
\begin{equation}\label{3.7}
\begin{array}{l}
\eta(|\xi|^2+\sigma^2+|\tau|^\alpha)^{-1/2}|\tilde{p}_\psi|^2+2\{\Re \tilde{p}_\psi, \Im \tilde{p}_\psi\}\gtrsim (|\xi|^2+\sigma^2+|\tau|^\alpha)^{3/2}.
\end{array}
\end{equation}
By G{\aa}rding's inequality and \eqref{3.6}, we obtain for large enough $\eta$ that
\begin{equation}\label{3.8}
\begin{array}{rl}
||\chi_\nu u||^2_{H^{3\alpha/4,3/2}}
\lesssim||P_\psi\chi_\nu u||_{L^2}^2+||u||^2_{H^{\alpha/2,1}}
\end{array}
\end{equation}
which proves \eqref{3.2}.

\section{Carleman estimates}\label{sec4}

In this section, we will derive a Carleman estimate from \eqref{3.1} by conjugating $u$ in \eqref{3.1} by $e^{i\beta z}$ with
the large parameter $\beta$.

\begin{lemma}\label{lem4.1}
There exist a sufficiently large constant $\beta_1$ depending on $n$ such that for all $v(t,x)\in C_0^{\infty}(U)\cap\dot{\mathcal{S}}(\bar{\mathbf{R}}_+^{1+n})$ and $\beta\geq \beta_1$, we have that
\begin{equation}\label{4.1}
\begin{array}{l}
\sum_{|\gamma|\leq1}\beta^{3-2|\gamma|}\int e^{2\beta\psi(x)}|D_x^\gamma v|^2dtdx
\lesssim \int e^{2\beta\psi(x)}|P(t,x;D_t,D_x) v|^2dtdx.
\end{array}
\end{equation}
\end{lemma}

\pf. Let $u\in C^\infty_0(U\times(-z_0,z_0))$ and $\beta\in \mathbf{R}$. We denote
$$u\hat{}(\sigma,t,x)=\int u(z,t,x)e^{-i\sigma z}dz.$$
Then
\begin{equation}\label{4.2}
\begin{array}{rl}
e^{-i\beta z}P_\psi (e^{i\beta z}u) =\int e^{i(\sigma-\beta)z}p(t,x,D_t,D_x+i|\sigma|\nabla\psi)u\hat{}(\sigma-\beta,t,x)d\sigma.
\end{array}
\end{equation}
By the Leibniz formula, it yields that for any smooth function $w=w(t,x)$,
\begin{equation}\label{4.3}
\begin{array}{rl}
&p(t,x,D_t,D_x+i|\sigma|\nabla\psi)w\\
=&e^{|\sigma|\psi}p(t,x,D_t,D_x) e^{-|\sigma|\psi}w\\
=&e^{(|\sigma|-|\beta|)\psi}e^{|\beta|\psi}p(t,x,D_t,D_x) e^{-(|\sigma|-|\beta|)\psi}e^{-|\beta|\psi}w\\
=&p(t,x,D_t,D_x+i|\beta|\nabla\psi)w+\Sigma_{k+|\gamma|<2,j>0}C_{j,k,\gamma}(x)(|\sigma|-|\beta|)^j\beta^kD_x^\gamma w\\
&+C_0(|\sigma|-|\beta|)D_t^{\alpha-1} w,
\end{array}
\end{equation}
where the last term $C_0(|\sigma|-|\beta|)D_t^{\alpha-1} w$ is coming from the last term of \eqref{2.2} whose symbol is $\frac{X}{T}i^\alpha(\tau+i\tau_0)^{\alpha-1}\xi_n$.

Now, let $u(z,t,x)=f(t,x)g(z)$, then we have $u\hat{}(\sigma-\beta,t,x)=f(t,x)\hat{g}(\sigma-\beta)$.
Applying \eqref{4.3} to $w=f(t,x)\hat{g}(\sigma-\beta)$, we have from \eqref{4.2} and \eqref{4.3}
\begin{equation}\label{4.4}
\begin{array}{rl}
&e^{-i\beta z}P_\psi (e^{i\beta z}f(t,x)g(z))\\
=&\int e^{i(\sigma-\beta)z}p(t,x,D_t,D_x+i|\sigma|\nabla\psi)f(t,x)\hat{g}(\sigma-\beta)d\sigma\\
=&\int e^{i(\sigma-\beta)z}p(t,x,D_t,D_x+i|\beta|\nabla\psi)f(t,x)\hat{g}(\sigma-\beta)d\sigma\\
&+\int e^{i(\sigma-\beta)z}\Sigma_{j+k+|\gamma|\le 2,\,j>0}C_{j,k,\gamma}(x)(|\sigma|-|\beta|)^j\beta^kD_x^\gamma f(t,x)\hat{g}(\sigma-\beta)d\sigma\\
&+\int e^{i(\sigma-\beta)z}C_0(|\sigma|-|\beta|)D_t^{\alpha-1} f(t,x)\hat{g}(\sigma-\beta)d\sigma\\
=&g(z)p(t,x,D_t,D_x+i|\beta|\nabla\psi)f(t,x)\\
&+\Sigma_{j+k+|\gamma|\le 2,\,j>0}C_{j,k,\gamma}(x)G_j(\beta)g(z)\beta^kD_x^\gamma f(t,x)\\
&+C_0D_t^{\alpha-1} f(t,x)G_1(\beta)g(z),
\end{array}
\end{equation}
where $G_j(\beta)g(z)=\int e^{i(\sigma-\beta)z}(|\sigma|-|\beta|)^j\hat{g}(\sigma-\beta)d\sigma
=\int e^{i\sigma z}(|\sigma+\beta|-|\beta|)^j\hat{g}(\sigma)d\sigma$.
By the Plancherel theorem, we have
\begin{equation}\label{4.5}
\begin{array}{l}
||G_j(\beta)g(z)||^2\lesssim \int |\sigma|^{2j}|\hat{g}(\sigma)|^2d\sigma\lesssim||g||^2_{H^j(\mathbf{R})}.
\end{array}
\end{equation}
Let $g\in C_0^\infty((-z_0,z_0))$ be any non-zero function. Combining \eqref{4.4} and \eqref{4.5}, it implies that
\begin{equation}\label{4.6}
\begin{array}{rl}
&||P_\psi (e^{i\beta z}f(t,x)g(z))||\\
\lesssim&||g||\cdot||p(t,x,D_t,D_x+i|\beta|\nabla\psi)f||+\Sigma_{j+k+|\gamma|\le 2,j>0}\beta^k\cdot||g||_{H^j(\mathbf{R})}\cdot||D_x^\gamma f||\\
&+||D_t^{\alpha-1} f||\cdot||g||_{H^1(\mathbf{R})}\\
\lesssim&||p(t,x,D_t,D_x+i|\beta|\nabla\psi)f||+\Sigma_{k+|\gamma|<2}\beta^k\cdot||D_x^\gamma f||+||D_t^{\alpha-1} f||.
\end{array}
\end{equation}

On the other hand, we need to estimate the lower bound on the left hand side of \eqref{3.1}.
We let $u(t,x,z)=e^{i\beta z}f(t,x)g(z)$ in \eqref{3.1}. A direct computation gives that
\begin{equation}\label{4.7}
\begin{array}{rl}
e^{-i\beta z}h(D_z)^jD_z^k(e^{i\beta z}g(z))&=(2\pi)^{-1}\int e^{i(\sigma-\beta)z}h(\sigma)^j\sigma^k\hat{g}(\sigma-\beta)d\sigma\\
&=(2\pi)^{-1}\int e^{i\sigma z}h(\sigma+\beta)^j(\sigma+\beta)^k\hat{g}(\sigma)d\sigma\\
&=(2\pi)^{-1}\int e^{i\sigma z}[h(\sigma+\beta)-h(\beta)+h(\beta)]^j(\sigma+\beta)^k\hat{g}(\sigma)d\sigma\\
&=h(\beta)^j\beta^kg(z)+\Sigma_{j',k'\in J}C_{j',k'}h(\beta)^{j-j'}\beta^{k-k'}H_{j',k'}(\beta)g(z),
\end{array}
\end{equation}
where $J=\{(j',k'):0\leq j'\leq j,0\leq k'\leq k\}\setminus(0,0)$ and
$H_{j',k'}(\beta)g(z)=\int e^{i\sigma z}[h(\sigma+\beta)-h(\beta)]^{j'}\sigma^{k'}\hat{g}(\sigma)d\sigma$.
By \eqref{4.7}, we have for large $\beta$ that
\begin{equation}\label{4.8}
\begin{array}{l}
\begin{cases}
||e^{-i\beta z}h(D_z)^jD_z^k(e^{i\beta z}g(z))-h(\beta)^j\beta^kg(z)||\lesssim h(\beta)^j\beta^{k}(h(\beta)^{-1}+\beta^{-1})\\
||h(D_z)^jD_z^k(e^{i\beta z}g(z))||\gtrsim h(\beta)^{j}\beta^{k}.
\end{cases}
\end{array}
\end{equation}
By \eqref{4.8}, we have for large $\beta$
\begin{equation}\label{4.9}
\begin{array}{rl}
&||h(D_z)^j\Lambda_\alpha^{s}  D_z^k(u)||\\
=&||\Lambda_{\alpha}^{s}  f(t,x)e^{-i\beta z}h(D_z)^jD_z^k(e^{i\beta z}g(z))||\\
=&||\Lambda_{\alpha}^{s}  f(t,x)[e^{-i\beta z}h(D_z)^jD_z^k(e^{i\beta z}g(z))-h(\beta)^j\beta^kg(z)+h(\beta)^j\beta^kg(z)]||\\
\gtrsim&h(\beta)^{j}\beta^{k}||\Lambda_{\alpha}^{s}  f||.
\end{array}
\end{equation}
Recall that $h(\beta)\simeq\beta^{1/2}$, by
\eqref{4.9} and \eqref{3.1}, we obtain
\begin{equation}\label{4.10}
\begin{array}{rl}
\sum_{|\gamma|\leq1}\beta^{3-2|\gamma|}\int |D_x^\gamma f|^2dtdx&\lesssim\Sigma_{k+s<2}h(\beta)^{2(2-k-s)}\beta^{2k}||\Lambda_{\alpha}^{s} f||^2\\
&\lesssim\Sigma_{k+s<2}||h(D_z)^{2-k-s}\Lambda_{\alpha}^{s} D_z^k(u)||^2\\
&\lesssim||P_\psi u||^2\\
&=||P_\psi(e^{i\beta z}f(t,x)g(z))||^2.
\end{array}
\end{equation}

Combining \eqref{4.10} and \eqref{4.6}, we have for large enough $\beta$ that
\begin{equation}\label{4.11}
\begin{array}{rl}
\sum_{|\gamma|\leq1}\beta^{3-2|\gamma|}\int |D_x^\gamma f|^2dtdx
&\lesssim||p(t,x,D_t,D_x+i|\beta|\nabla\psi)f||^2.
\end{array}
\end{equation}
By letting $f=e^{\beta \psi}v$ in \eqref{4.11}, we immediately have \eqref{4.1}.
\eproof

\section{Proof of Theorem \ref{thm1.1}}\label{sec5}

This section is devoted to the proof of the main theorem, Theorem \ref{thm1.1}.
Since $u\in H^{\alpha,2}(\mathbf{R}^{1+n})$ and $u(t,y)=0\quad (t\leq0)$, we can find a sequence $\{u_m\}$ in $\mathcal{S}(\bar{\mathbf{R}}_+^{1+n})$
which converges to $u$ in $H^{\alpha,2}(\mathbf{R}^{1+n})$. The limit arguments in \eqref{4.1} imply that we can assume $u\in \mathcal{S}(\bar{\mathbf{R}}_+^{1+n})$.
To apply Lemma \ref{lem4.1}, we first recall
$${\rm supp}\, u\subset\{x_n\geq -X\}$$
and then define a smooth function $\chi$ by

\begin{equation}\label{5.1}
\chi (x_n)=
\begin{cases}
\begin{array}{l}
1,\quad x_n\leq X/2,\\
0,\quad x_n\geq X.
\end{array}
\end{cases}
\end{equation}
From \eqref{2.1}, it is not hard to get $\chi u\in C_0^{\infty}(U)\cap\dot{\mathcal{S}}(\bar{\mathbf{R}}_+^{1+n})$.
Thus, we can apply $\chi u$ to the Carleman estimates \eqref{4.1} and get that
\begin{equation}\label{5.2}
\begin{array}{rl}
&\sum_{|\alpha|\leq1}\beta^{3-2|\gamma|}\int_{x_n\leq X/2} e^{2\beta\psi(x)}|D^\alpha u|^2dtdx\\
\leq &\sum_{|\alpha|\leq1}\beta^{3-2|\gamma|}\int e^{2\beta\psi(x)}|D^\gamma (\chi u)|^2dtdx\\
\lesssim &\int e^{2\beta\psi(x)}|P(t,x;D_t,D_x) (\chi u)|^2dtdx\\
\lesssim &\sum_{|\alpha|\leq1}\int_{x_n\leq X/2} e^{2\beta\psi(x)}|D^\gamma u|^2dtdx+\int_{X/2<x_n\leq X} e^{2\beta\psi(x)}|[P,\chi] u|^2dtdx,
\end{array}
\end{equation}
where $[\cdot,\cdot]$ denotes the commutator. Let $\beta$ be large enough to absorb the first term on the right hand side of \eqref{5.2}, we get from \eqref{5.2} that
\begin{equation}\label{5.3}
\begin{array}{rl}
&\beta^{3}\int_{x_n\leq X/4} e^{9\beta X^2/16}|u|^2dtdx\\
\leq &\sum_{|\alpha|\leq1}\beta^{3-2|\alpha|}\int_{x_n\leq X/2} e^{2\beta\psi(x)}|D^\alpha u|^2dtdx\\
\lesssim &C(u)e^{\beta X^2/4}.
\end{array}
\end{equation}
Let $\beta$ tend to $\infty$, we obtain that $u=0$ on $x_n\leq X/4$.
Recall that $x_n=y_n+|y'-\hat{y}'|^2+\frac{X}{T}(t-T)$, the proof is complete.

\end{document}